\documentstyle{article}
\newtheorem{theorem}{Theorem}
\newtheorem{prop}[theorem]{Proposition}
\newtheorem{cor}[theorem]{Corollary}
\newtheorem{example}[theorem]{Example}
\newtheorem{definition}[theorem]{Definition}
\newtheorem{remark}[theorem]{Remark}
\title{Singular perturbations, regularization and extension theory}
\author{H. Neidhardt \and V.A.Zagrebnov} \date{}

\begin{document}
\maketitle
For nonpositive singular potentials in quantum mechanics it can
happen that the Schr{\"o}dinger operator is not essentially
self-adjoint on a natural domain of definition or not semibounded
{}from below. In this case we have a lot of self-adjoint extensions
each of them is a candidate for the right physical Hamiltonian of
the system. Hence the problem arises to single out the right
physical self-adjoint extension. Usually this problem is solved as
follows. At first one has to approximate the singular potential by
a sequence of bounded potentials (cut-off approximation). After
that one has to show that the arising sequence of Schr\"odinger
operators converges in the strong resolvent sense to one of the
self-adjoint extensions if the cut-off approximation tends to the
singular potential. The so determined self-adjoint extensions is
regarded as the right physical Hamiltoninan. Very often the right
physical Hamiltonian coincides with the Friedrichs extension.

With respect to the Schr{\"o}dinger operator in $L^2({\bf R}^2)$
this problem was discussed by \cite{G}, \cite{H}, \cite{K},
\cite{Sch} and \cite{S}. An operator-theoretical investigation of
this problem was started by Nenciu in \cite{N} and continued by the
authors in \cite{NZ}. In the following we continue those abstract
investigations. We assume that a semibounded symmetric operator
admits a monotonously decreasing sequence of semibounded symmetric
operators such that the corresponding sequence of Friedrichs
extensions converges in the strong resolvent sense to the
Friedrichs extension of the symmetric operator with which we have
started. The problem will be to find necessary and sufficient
conditions that any other sequence of semibounded self-adjoint
extensions of the decreasing sequence of symmetric operators
converges to this Friedrichs extension too. Unfortunately, we are
unable to solve ths problem in full generality. This means we have
found a necessary condition which must be satisfied in order to
have the desired convergence. However, we can prove the converse
only for special sequences of self-adjoint extensions but not for
all.

In more detail the problem can be described as follows. Let $A$ and
$V$ be two nonnegative self-adjoint operators on the separable
Hilbert space ${\cal H}$. Further, let ${\cal D} \subseteq dom(A)
\cap dom(V)$ a dense subset of ${\cal H}$ such that
\begin{equation}
(Vf,f) \le a(Af,f) + b\|f\|^2, \qquad f \in {\cal D}, \qquad 0 <
a,b. \end{equation}
We introduce the abstract operator $\dot{H}_{\alpha}$
\begin{equation}
\dot{H}_{\alpha}f = Af - {\alpha}Vf, \quad f \in
\mbox{dom}(\dot{H}_{\alpha}) = {\cal D}, \quad {\alpha} > 0.
\end{equation}
If the coupling constant ${\alpha}$, ${\alpha > 0}$, obeys
${\alpha} < 1/a$, then the operator $\dot{H}_{\alpha}$ is
symmetric, closable and semibounded with lower bound $-{\alpha}b$.
However, the operator $\dot{H}_{\alpha}$ is in general not
esssentially self-adjoint.
\begin{example}\label{A}
{\em Let ${\cal H} = L^2({\bf R}^1)$ and let $A$ be the usual
Laplace operator on $L^2({\bf R}^1)$, i.e. $A = -d^2/dx^2$. By $V$
we denote the multiplication operator arising from the real
potential $V(x)$,
\begin{equation}
V(x) = \frac{1}{4\kappa}\frac{1}{|x|^{\beta}}, \qquad 1 \le \beta
\le 2, \qquad \kappa > 0. \end{equation}
Let ${\cal D} = C^\infty_0({\bf R}^1\setminus\{0\})$. If $1 \le
\beta < 2$, then for every $\kappa > 0$ there are real numbers $a <
1$ and $b \ge 0$ such that
\begin{equation}
\int^\infty_{-\infty}
\frac{1}{4\kappa}\frac{1}{|x|^{\beta}}|f(x)|^2dx \le
a\int^\infty_{-\infty}|f'(x)|^2dx + b
\int^\infty_{-\infty}|f(x)|^2dx. \end{equation}
for $\kappa > 0$. If $\beta = 2$, then this is only true for
$\kappa > 1$.}
\end{example}
\begin{example}\label{B}
{\em Let ${\cal H} = L^2({\bf R}^2)$ and let $A$ be the usual
Laplace operator on $L^2({\bf R}^2)$, i.e. $A = -\Delta$. Further,
let $\Gamma$ be a smooth curve in ${\bf R}^2$ which is
parameterized by
\begin{equation}
\Gamma = \{(x,y) \in {\bf R}^2: x = \rho(\varphi)\cos\varphi, y =
\rho(\varphi)\sin\varphi, 0 \le \varphi < 2\pi\} \end{equation}
where $\rho(\varphi) > 0$ is a smooth function. Again $V$ is the
multiplication operator arising from
\begin{equation}
V(x) = \frac{1}{5\kappa}\frac{1}{{||x| - \rho(\varphi)|}^{\beta}},
\quad \quad 1 \le \beta \le 2, \quad |x| = \sqrt{x^2 + y^2}.
\end{equation}
We set ${\cal D} = C^\infty_0({\bf R}^2\setminus \Gamma)$. If $1
\le \beta < 2$, then for every $\kappa > 0$ there are real numbers
$a < 1$ and $b \ge 0$ such that
\begin{equation}
\int_{{\bf R}^2} \frac{1}{5\kappa}\frac{1}{{||x| -
\rho(\varphi)|}^{\beta}}|f(x)|^2dx \le a\int_{{\bf R}^2}|\nabla
f(x)|^2dx + b\int_{{\bf R}^2}|f(x)|^2dx. \end{equation}
For $\beta = 2$ this is true only for $\kappa > 1$.} \end{example}
Let us assume that the $\dot{H}_{\alpha}$ is not essentially
self-adjoint. Since $\dot{H}_{\alpha}$ is semibounded the
Friedrichs extension $\hat{H}_{\alpha}$ exists. Moreover, denoting
by $\hat{A}$ the Friedrichs extension of $\dot{A} = A|{\cal D}$ it
is not hard to see that $\hat{H}_{\alpha}$ coincides with the form
sum of $\hat{A}$ and $-{\alpha}V$, i.e.
\begin{equation}
\hat{H}_{\alpha} = \hat{A} \dot{+} (-{\alpha}V). \end{equation}
In the above examples the Friedrichs extension corresponds to the
Dirichlet boundary condition at $x = 0$ for the first example and
on $\Gamma$ for the second one.

Next let us introduce a regularizing sequence for the singular
perturbation.
\begin{definition}
{\em A sequence $\{V_n\}_{n \ge 1}$ of bounded non-negative
self-adjoint operators is called a regularizing sequence of $V$ if}
\begin{enumerate}
\item[(i)] $V_1 \le V_2 \le \ldots \le V_n \le \dots \le V$
\item[(ii)] $lim_{n\to\infty}(V_nf,f) = (Vf,f)$, $f \in {\cal D}
\subseteq dom(V)$. \end{enumerate}
\end{definition}
\begin{example}\label{C}
{\em In the Examples \ref{A} and \ref{B} the sequence $V_n$ is
given as multiplication operators with the cut-off potentials}
\begin{equation}
V_n(x) = \inf_{x \in {\bf R}^l} \{n,V(x)\}, \qquad l =1,2.
\end{equation}
\end{example}
With the regularizing sequence $\{V_n\}^{\infty}_{n=1}$ we
associate the following sequence of self-adjoint operators
$H_{\alpha,n}$,
\begin{equation}
H_{\alpha,n} = A - {\alpha}V_n, \qquad n = 1,2,\ldots.
\end{equation}
The problem is now to find conditions which guarantee that the
approximating sequence $\{H_{\alpha,n}\}^{\infty}_{n=1}$ tends to
the Friedrichs extension $\hat{H}_{\alpha}$, i.e.,
\begin{equation}
s-{\lim_{n\to\infty}}(H_{\alpha,n} - z)^{-1} = (\hat{H}_{\alpha} -
z)^{-1}, \qquad \Im\mbox{m}(z) \not= 0 \end{equation}

However, from the mathematical point of view this setup seems to be
unnatural. To explain this we remark that for any $n = 1,2,\ldots$
the operator $H_{\alpha,n}$ is a self-adjoint extension of the
semibounded symmetric operator $\dot{H}_{\alpha,n} =
H_{\alpha,n}|{\cal D} = \dot{A} - {\alpha}V_n$, i.e.
$\dot{H}_{\alpha,n} \subseteq H_{\alpha,n}$. Taking another
semibounded self-adjoint extension $\tilde{A}$ of $\dot{A}$ we get
another sequence $\tilde{H}_{\alpha,n}$,
\begin{equation}
\tilde{H}_{\alpha,n} = \tilde{A} - {\alpha}V_n, \qquad n
=1,2,\ldots, \end{equation}
which naturally implies the question: why we should to investigate
the convergence for $H_{\alpha,n}$ and why not for
$\tilde{H}_{\alpha,n}$? So in the following we shall search for
conditions which guarantee that
\begin{equation}\label{conv}
s-{\lim_{n\to\infty}}(\tilde{H}_{\alpha,n} - z)^{-1} =
(\hat{H}_{\alpha} - z)^{-1}, \qquad \Im\mbox{m}(z) \not= 0.
\end{equation}
for any semibounded self-adjoint extension $\tilde{A}$ of
$\dot{A}$. In particular, this would be clarified the uniqueness
problem of the limit (\ref{conv}) for the two ''extreme cases'':
the sequence of Friedrichs extension $\hat{H}_{\alpha,n}$,
\begin{equation}
\hat{H}_{\alpha,n} = \hat{A} - {\alpha}V_n, \qquad n = 1,2,\ldots,
\end{equation}
where $\hat{A}$ is the Friedrichs extension of $\dot{A}$, and of
the sequence of Krein extensions $\check{H}_{\alpha,n}$
\begin{equation}
\check{H}_{\alpha,n} = \check{A} - {\alpha}V_n, \qquad n =
1,2,\ldots, \end{equation}
where $\check{A}$ is the Krein extension (soft extension)
\cite{AS}, \cite{BN}, \cite{Kr} of $\dot{A}$ with respect to a
given lower bound $\eta < 0$, i.e. $\check{A} \ge {\eta}I$.

In general we cannot expect that the sequence
$\tilde{H}_{\alpha,n}$ tends to $\hat{H}_{\alpha}$ assuming only
that $\{V_n\}_{n\ge1}$ is a regularizing sequence. Actually we need
a little bit more. Only if $\tilde{A}$ is the Friedrichs extension
$\hat{A}$ of $\dot{A}$, i.e. $\tilde{A} = \hat{A}$, then we obtain
\begin{equation}
s-{\lim_{n\to\infty}}(\hat{H}_{\alpha,n} - z)^{-1} =
(\hat{H}_{\alpha} - z)^{-1}, \qquad \Im\mbox{m}(z) \not= 0,
\end{equation}
without any additional assumptions \cite{NZ}. How to find this
additional assumptions? An essential hint comes from the following
proposition.
\begin{prop}\label{necessprop}
Let $\{V_n\}_{n\ge1}$ be a regularizing sequence of $V$. If for
every self-adjoint extension $\tilde{A}$ of $\dot{A} = A|{\cal D}$
obeying $\tilde{A} \ge \eta$, $\eta < 0$, the convergence
(\ref{conv}) takes place, then
\begin{equation}\label{necess}
\sup_{n \ge 1}(V_{n}h,h) = +\infty
\end{equation}
for every nontrivial $h$ of ${\cal N}_{\eta} = ker(\dot{A}^\ast -
\eta)$. \end{prop}
By this proposition it seems to be natural to introduce the
following notation.
\begin{definition}
{\em Let $\{V_n\}_{n \ge 1}$ be a regularizing sequence of $V$. The
sequence is called admissible with respect to $\dot{A} = A|{\cal
D}$ if there is a $\eta < 0$ such that for every nontrivial $h \in
{\cal N}_{\eta} = ker(\dot{A}^\ast - \eta)$ the condition
(\ref{necess}) is satisfied.}
\end{definition}
\begin{remark}
{\em It can be shown that if (\ref{necess}) is satisfied for one
$\eta< 0$, then it holds for every $\eta' < 0$. So the property
(\ref{necess}) is independent on $\eta < 0$.} \end{remark}
\begin{example}
{\em It can be shown that the regularizing sequences of Example
\ref{C} for the Examples \ref{A} and \ref{B} are admisssible with
respect to $\dot{A} = -\frac{d}{dx^2}|C^{\infty}_{0}({\bf R}^1
\setminus \{0\})$ and $\dot{A} = -\Delta|C^{\infty}_{0}({\bf R}^2
\setminus \Gamma)$.} \end{example}
Hence, the optimal way to solve our problem would be to show that
the converse to Proposition \ref{necessprop} is true, i.e., if
$\{V_n\}_{n\ge 1}$ is an admissible regularizing sequence of $V$
with respect to $\dot{A} = A|{\cal D}$, then for every semibounded
self-adjoint extension $\tilde{A}$ of $\dot{A}$ we have that the
convergence (\ref{conv}) is valid. Till now we cannot prove this
conjecture in full generality. However, if we restrict the set of
semibounded self-adjoint extensions $\tilde{A}$ of $\dot{A}$, then
we can do it. To describe these restrictions we use a description
of all semibounded self-adjoint extensions which goes back to
\cite{AS}. Let $\tilde{A}$ be any semibounded self-adjoint
extension of $\dot{A} = A|{\cal D}$ with lower bound greater than
$\eta < 0$, i.e. $\tilde{A} \ge \eta$. By $\tilde{\nu} \ge \eta$ we
denote the closed quadratic form which corresponds to $\tilde{A}$,
i.e.
\begin{eqnarray}\label{corresponding form} \tilde{\nu}(f,f) & = &
(({\tilde{A} - \eta})^{1/2}f,({\tilde{A} - \eta})^{1/2}f) +
\eta(f,f),\\ f \in \mbox{dom}(\tilde{\nu}) & = &
\mbox{dom}(({\tilde{A} - \eta})^{1/2}). \nonumber \end{eqnarray}
In particular, by $\hat{\nu} \ge 0$ we denote the closed quadratic
form which corresponds to the Friedrichs extension $\hat{A}$ of
$\dot{A}$. In accordance with \cite{AS} we have an one-to-one
correspondence between the set of all semibounded self-adjoint
extensions $\tilde{A}$ of $\dot{A}$ obeying $\tilde{A} \ge \eta$
and all non-negative closed quadratic forms $\tilde{q}$ on the
deficiency subspace ${\cal N}_{\eta} = ker(\dot{A}^\ast - \eta)$,
where the form $\tilde{q}$ is not necessarily densely defined on
${\cal N}_{\eta}$. The correspondence is given by the formulas
\begin{equation}\label{domain}
\mbox{dom}(\tilde{\nu}) = \mbox{dom}(\hat{\nu}) \dot{+}
\mbox{dom}(\tilde{q}), \end{equation}
where $\dot{+}$ means $\mbox{dom}(\hat{\nu}) \cap
\mbox{dom}(\tilde{q}) = \{0\}$, and
\begin{equation}\label{form}
\tilde{\nu}(g+h,g+h) = \hat{\nu}(g,g) + \tilde{q}(h,h) +
2\eta\Re(g,h) + \eta(h,h), \end{equation}
$g \in \mbox{dom}(\hat{\nu}), h \in \mbox{dom}(\tilde{q}) \subseteq
{\cal N}_{\eta}$. Therefore, starting with extension $\tilde{A}$
which obeys $\tilde{A} \ge \eta$ we can find a unique non-negative
closed quadratic form $\tilde{q}$ on ${\cal N}_{\eta}$ such that
(\ref{domain}) and (\ref{form}) holds. Conversely, if we have a
non-negative closed quadratic from $\tilde{q}$ on ${\cal
N}_{\eta}$, then we can define by (\ref{domain}) and (\ref{form}) a
semibounded extension $\tilde{A}$ of $\dot{A}$ obeying $\tilde{A}
\ge \eta$. The domain of $\tilde{q}$ may be a closed subspace of
${\cal N}_{\eta}$ or not. The Friedrichs extension $\hat{A}$
corresponds to the trivial form $\hat{q}$, i.e.,
$\mbox{dom}(\hat{q}) = \{0\}$. Very often this is expressed by
$\hat{q} = +\infty$. The Krein extension (soft extension)
\cite{AS}, \cite{BN}, \cite{Kr} $\check{A}$ with respect to the
lower bound $\eta < 0$, i.e. $\check{A} \ge {\eta}I$, is given by
the form $\check{q}$ which is zero on the whole deficiency subspace
${\cal N}_{\eta}$, i.e., $\check{q} = 0$. All other forms
$\tilde{\nu}$ are between $\check{\nu}$ and $\hat{\nu}$ which
yields $\check{A} \le \tilde{A} \le \hat{A}$.

Of course the description is only unique if we fix some $\eta < 0$.
Changing $\eta$ we get different quadratic forms $\tilde{q}_{\eta}$
for the same semibounded self-adjoint extension $\tilde{A}$ of
$\dot{A}$. However, there are some invariants which do not depent
on $\eta$. For instance, if $\mbox{dom}(\tilde{q}_{\eta})$ is a
closed subspace in ${\cal N}_{\eta}$, then
$\mbox{dom}(\tilde{q}_{\eta'})$ is a closed subspace for $\eta'(<
0)$, too.

Using this description our main theorem can be formulated now as
follows.
\begin{theorem}
Let $\{V_n\}_{n \ge 1}$ be an admissible regularizing sequence of
$V$ with respect to $\dot{A}$ and let $\tilde{A}$ be a self-adjoint
extension of $\dot{A}$ obeying $\tilde{A} \ge \eta$ for some $\eta
< 0$. If $\tilde{A}$ corresponds to a closed quadratic form
$\tilde{q}$ on ${\cal N}_{\eta} = \mbox{ker}({\dot{A}}^{\ast} -
\eta)$ and the domain $\mbox{dom}(\tilde{q})$ is a closed subspace
of ${\cal N}_{\eta}$, then for sufficiently small coupling
constants ${\alpha} > 0$ we have
\begin{equation}\label{convergence}
s-\lim_{n\to\infty} (\tilde{H}_{\alpha,n} - z)^{-1} =
(\hat{H}_{\alpha} - z)^{-1}, \qquad \Im\mbox{m}(z) \not= 0,
\end{equation}
where $\hat{H}_{\alpha}$ is the Friedrichs extension of
$\dot{H}_{\alpha} = (A - {\alpha}V)|{\cal D}$.

In particular, if $\check{A}$ denotes the Krein extension of
$\dot{A}$ with respect to the lower bound $\eta < 0$, then for
sufficiently small $\alpha > 0$ we have
\begin{equation}
s-\lim_{n\to\infty} (\check{H}_{\alpha,n} - z)^{-1} =
(\hat{H}_{\alpha} - z)^{-1}, \qquad \Im\mbox{m}(z) \not= 0.
\end{equation}
\end{theorem}
If the deficieny indices are finite, then the theorem admits a
strengthening.
\begin{theorem}\label{finite}
If the deficieny indices of $\dot{A}$ are finite, then for any
self-adjoint extension $\tilde{A}$ of $\dot{A}$ and any coupling
constant ${\alpha} < 1/a$ we have (\ref{convergence}).
\end{theorem}
The Theorem \ref{finite} improves the results of Section 3 of
\cite{NZ}. Moreover, the theorem can be slightly generalized.
\begin{cor}\label{pseudofinite}
If $\tilde{A}$ is a semibounded self-adjoint extension of $\dot{A}$
such that
\begin{equation}
\mbox{dim}(\mbox{dom}(\tilde{\nu}) / \mbox{dom}(\hat{\nu})) <
+\infty, \end{equation}
then for ${\alpha} < 1/a$ (\ref{convergence}) is valid. \end{cor}
The theorems and corollary admit an application to our examples.
\begin{example}
{\em Since in Example \ref{A} the deficiency indices of $\dot{A} =
-\frac{d}{dx^2}|C^{\infty}_{0}({\bf R}^1 \setminus \{0\})$ are
finite by Corollary \ref{finite} we always have the desired
convergence (\ref{convergence}).

In Example \ref{B} we have the desired convergence
(\ref{convergence}) only for a special set of self-adjoint
extensions of $\dot{A} = -\Delta|C^{\infty}_{0}({\bf R}^2 \setminus
\Gamma)$. The set includes the Krein extension (the corresponding
boundary condition can be found in \cite{AS}) and extensions which
are characterized by Corollary \ref{pseudofinite}. However, it
remains an open question: whether the sequence of usual
Schr\"odinger operatos $H_{\alpha,n} = -\Delta - {\alpha}V_n$,
where $-\Delta$ denotes the usual Laplace operator in $L^2({\bf
R}^2)$ convergences to the Friedrichs extension of the symmetric
operator $(-\Delta - {\alpha}V)|C^{\infty}_{0}({\bf R}^2)$? The
problem is that the domain of the closed quadratic form, which by
(\ref{corresponding form}) - (\ref{form}) corresponds to the usual
Laplace operator $-\Delta$ in $L^2({\bf R}^2)$ regarded as a
self-adjoint extension of $-\Delta|C^{\infty}_{0}({\bf R}^2)$, is
not a closed subspace in ${\cal N}_{\eta}$.} \end{example}
\begin{remark}
{\em If the deficiency indices are finite, then the strong
resolvent convergence (\ref{convergence}) can be replaced by the
operator-norm convergence \cite{NZ}. However, if the deficiency
indices are infinite this is not true in general. For instance, let
in Example \ref{B} the curve $\Gamma$ be the unite circle. Then one
can show that for any interval $\delta \subseteq (-\infty,0)$ and
any integer $N$ there is a greater integer $n \ge N$ such that
$H_{\alpha,n}$ has an eigenvalue in $\delta$. Consequently, this
excludes the operator-norm convergence for the operators
$\{H_{\alpha,n}\}_{n \ge 1}$.} \end{remark}
\end{document}